\documentclass[12pt, fleqn]{amsart}

\DeclareMathOperator{\cov}{cov}
\DeclareMathOperator{\Var}{Var}
\DeclareMathOperator{\cor}{corr}

\usepackage[square,compress,comma, numbers,sort]{natbib}
\usepackage[colorlinks=true, citecolor=red, linkcolor=blue]{hyperref}
\usepackage{amsfonts,mathtools}
\usepackage{amsmath}

\allowdisplaybreaks[4]

\def\ve{\varepsilon}

\usepackage{amssymb}
\usepackage{color}
\usepackage{float}

\newcounter{alphasect}
\def\alphainsection{0}

\let\oldsection=\section
\def\section{%
  \ifnum\alphainsection=1%
    \addtocounter{alphasect}{1}
  \fi%
\oldsection}%

\renewcommand\thesection{%
 \ifnum\alphainsection=1%
   \Alph{alphasect}%
 \else
   \arabic{section}%
 \fi%
}%

\usepackage{xparse}

\NewDocumentCommand{\ceil}{s O{} m}{%
	\IfBooleanTF{#1} 
	{\left\lceil#3\right\rceil} 
	{#2\lceil#3#2\rceil} 
}
\NewDocumentCommand{\floor}{s O{} m}{%
	\IfBooleanTF{#1} 
	{\left\lfloor#3\right\rfloor}
	{#2\lfloor#3#2\rfloor}
}

\numberwithin{equation}{section}
\newtheorem{theo}{Theorem}[section]
\newtheorem{sat}[theo]{Proposition}
\newtheorem{de}[theo]{Definition}
\newtheorem{lem}[theo]{Lemma}

\newtheorem{example}[theo]{Example}
\newtheorem{korr}[theo]{Corollary}
\newtheorem{remark}[theo]{Remark}

\numberwithin{equation}{section}

\newtheorem{lemma}{Lemma}[section]

\newcommand{\QED}{\hfill $\Box$}

\newcommand{\COM}[1]{}

\def\IF{\infty}

\newcommand{\R}{\mathbb{R}}
\newcommand{\inr}{\in \R}

\newcommand{\BQN}{\begin{eqnarray}}
\newcommand{\EQN}{\end{eqnarray}}
\newcommand{\BQNY}{\begin{eqnarray*}}
	\newcommand{\EQNY}{\end{eqnarray*}}

\newcommand{\limit}[1]{\lim_{#1 \to   \infty}}

\def\bqny#1{{\begin{eqnarray*} #1 \end{eqnarray*}}}
\def\bqn#1{{\begin{eqnarray} #1 \end{eqnarray}}}

\newcommand{\BS}{\begin{sat}}
	\newcommand{\ES}{\end{sat}}
\newcommand{\BT}{\begin{theo}}
	\newcommand{\ET}{\end{theo}}
\newcommand{\BK}{\begin{korr}}
	\newcommand{\EK}{\end{korr}}

\newcommand{\BEX}{\begin{example}}
	\newcommand{\EEX}{\end{example}}

\newcommand{\BD}{\begin{de}}
	\newcommand{\ED}{\end{de}}

\newcommand{\BIT}{\begin{itemize}}
	\newcommand{\EIT}{\end{itemize}}
\newcommand{\BDI}{\begin{description}}
	\newcommand{\EDI}{\end{description}}

\newcommand{\BRM}{\begin{remark}}
	\newcommand{\ERM}{\end{remark}}

\newcommand{\BEL}{\begin{lem}}
	\newcommand{\EEL}{\end{lem}}

\def\Z{\mathbb{Z}}

\def\fr{\pi_{\mathcal Z,H}(u)}

\newcommand{\pk}[1]{\mathbb{P} \left(#1 \right) }
\newcommand{\E}[1]{\mathbb{E} \left(#1 \right) }
\begin{document}

\title{Fractional Brownian Motion ruin model with random inspection time}

\author{{Grigori Jasnovidov}}
\address{Grigori Jasnovidov,
St. Petersburg Department
of Steklov Mathematical Institute
of Russian Academy of Sciences, St. Petersburg, Russia
}
\email{griga1995@yandex.ru}

 \maketitle

{\bf Abstract:} In this contribution we study the asymptotics of
\bqny{
\pk{\exists t\ge 0 : B_H(L(t))-cL(t)>u}, \quad u \to \IF,
}
where $B_H, H\in (0,1)$ is a fractional Brownian motion,
$L(t)$ is a non-negative pure jumps
L\'{e}vy process independent of $B_H$, $c>0$.



\section{Introduction}

The classical ruin probability
$$\pk{\exists t\ge 0: B(t)-ct>u}$$
with $B$ a standard Brownian motion, $u>0$ and $c>0$,
is a benchmark model in risk theory, see \cite{ig:69}.
This probability can be calculated explicitly, namely
it equals to $e^{-2cu}$ for all $c,u>0$, see \cite{DeM15}.
However, this is an exceptional case in the
field of ruin probabilities for Gaussian and related processes.
Typically probabilities of such type cannot be explicitly computed
and their asymptotics as $u \to \IF$ is dealt with, see, e.g.,
\cite{20lectures,PicandsA,DIEKER2005,Pit96}
and references therein.
One of possible extensions of the classical ruin probability is
$$\pk{\exists t\ge 0: B_H(t)-ct>u}$$
with $B_H, H \in (0,1)$ a standard fractional Brownian motion (fBm), i.e.,
a Gaussian process with zero mean and
$$\cov(B_H(t),B_H(s))=\frac{1}{2}(|t|^{2H}+|s|^{2H}-|t-s|^{2H}), \ \ \ \
 \ \ \ s,t \inr.$$
The ruin probability for fBm model and its extensions are a subject of numerous investigations,
see, e.g., \cite{HP2004,HP99,HA2013,Pit96,
ParisRuinGenrealHashorvaJiDebicki,KEP2015,HJ2014,PengLith}.
In this paper we consider the following
generalization of fBm ruin probability:
\bqn{\label{fBM_model}
\pi_{\mathcal Z,H}(u) = \pk{\exists t\ge 0 : B_H(L(t))-cL(t)>u},\quad u>0,
}
where $L(t)$ is a pure jumps L\'{e}vy process independent of $B_H$, $\mathcal Z$ is the non-negative distribution of jumps of $L$ and $c>0$.
Observe, that this model is an extension of discrete-time ruin probabilities
introduced in \cite{secondproj,Piterbargdiscrete}. Indeed,
taking $\mathcal Z$ to be a positive constant with probability one we reduce our model
to discrete ruin problems.
For studies devoted to L\'{e}vy time subordinations 
in the classical case $H=1/2$ we refer to \cite{levy_H=1/2_subordinated} and references therein.
We can rewrite \eqref{fBM_model} as follows:
\bqny{
\pi_{H,\mathcal Z}(u) = \pk{\exists t\ge 0:
B_H(t)-ct>u, \ t= Z_1+Z_2+...+Z_n \text{ for some } n},
}
where $Z_i$'s are i.i.d. with distribution $\mathcal Z$.
From the representation above it is easy to see that $L$ does not affect on $\pi_{H,\mathcal Z}(u)$, but $\mathcal Z$ may play a role. This additional randomness (compared to the classical case) brings supplementary complexities to our analysis. 
It turns out that cases $H<1/2$, $H=1/2$ and $H>1/2$
require different approaches and we treat them separately. For case $H>1/2$ the asymptotics of $\fr$ given in 
Theorem \ref{infinite-horizon_theorem} is as in the classical continuous case,
see \cite{HP99}, and $\mathcal Z$ does not effect on the answer.
This behavior is typical for long-range dependence case,
see also \cite{secondproj,ParisRuinGenrealHashorvaJiDebicki,
KEP20151,HJ2014,PengLith,KrzysPeng2015}.
For $H\le1/2$ we need to control behavior of $L(t)$
as $t$ tends to $\IF$, and thus we impose finiteness of variance (and hence expectation) of $\mathcal Z$:\\

\textbf{A}: $\mathbb{E}Z_1=\mu \in (0,\IF), \quad \Var(Z_1)<\IF$.\\

The asymptotics in Theorem \ref{infinite-horizon_theorem} for $H=1/2$ reminds the results in \cite{HashorvaGrigori}; 
the only difference is in Pickands-type constants.
For case $H<1/2$ due to short-range dependence of fBm
asymtpotics of $\fr$ as $u \to \IF$ coincides with the
asymptotics of the classical discrete ruin probability, see
\cite{secondproj,Debicki_Grisha}. The derivation of the 
results for case $H\le 1/2$ relies on seminal double-sum method, see, e.g., \cite{20lectures}.
\\
We organize the rest of the paper as follows. In the next 
section we present necessary tools and notation. Main
results of the paper are given in Section \ref{sec_main},
while the their proofs are relegated to Section \ref{sec_proofs}.

\section{Preliminaries}
To study the asymptotics of $\pi_{\mathcal Z,H}(u)$
it is worth to start with the asymptotics of
the classical continuous and discrete ruin probabilities. Let
$G(\delta) = \{0,\delta,2\delta,...\}$ for $\delta>0$ and $G(0) = [0,\IF)$.
Define for $u,c>0$ and $\delta\ge 0$
\bqny{
\psi_{\delta,H}(u)  = \pk{\exists t\in G(\delta) : B_H(t)-ct>u}
}
and $\psi_H(u) = \psi_{0,H}(u)$.
Define for $\delta\ge 0$ and $H\in (0,1)$ the
classical and discrete Pickands constants by
\bqny{
\mathbb{H}_{2H}^{\delta} = \limit{S}\frac{1}{S}
\E{\sup\limits_{t \in [0,S]\cap\delta\Z}e^{\sqrt 2B_H(t)-t^{2H}}}
}
and $\mathbb{H}_{2H} = \mathbb{H}_{2H}^{0}$.
It is known that $\mathbb{H}_{2H}^{\delta}\in (0,\IF)$
for all $\delta\ge0$ and $H\in (0,1)$, see
\cite{DiekerY,PicandsA,pickands_grisha}.
Let in the following
$$C_H = \frac{c^H}{H^H(1-H)^{1-H}}.$$
The proposition below establishes the asymptotic
behavior of $\psi_{\delta,H}(u)$ for
any $\delta\ge 0$:
\begin{sat}\label{proposition}
For any $\delta \ge 0$ as $u \to \IF$
\bqn{\label{one-dimensional_theorem} \ \ \ \ \ \ \ \
\psi_{\delta,H}(u) \sim
\begin{cases}
\mathbb{H}_{2H}\frac{2^{\frac{1}{2}-\frac{1}{2H}}\sqrt \pi}
{H^{1/2}(1-H)^{1/2}}
(C_Hu^{1-H})^{1/H-1}\Psi(C_Hu^{1-H}), & H>1/2 \text{ or } \delta=0,\\
\mathbb{H}^{2c^2\delta}_1  e^{-2cu}, & H=1/2 \text{ and } \delta>0\\
\frac{\sqrt{2\pi} H^{H+1/2} u^H}{\delta c^{H+1}(1-H)^{H+1/2}}
\Psi(C_Hu^{1-H}), & H<1/2 \text{ and } \delta>0.
\end{cases}
}
\end{sat}
The statement of the proposition above is given in Corollary 3.2 in \cite{ParisRuinGenrealHashorvaJiDebicki} (case $\delta =0$) and in Theorem 2.2 in \cite{secondproj}
(case $\delta>0$).

\section{Main Results}\label{sec_main}

Define for any positive distribution $\mathcal X$ the following
Pickands-type constant:
\bqny{
\mathcal H_{\mathcal X} =\limit{S}\frac{1}{S}
\int\limits_{\R}\pk{\exists t\ge 0:
 L_{\mathcal X}(t)\in [0,S],\sqrt 2
B( L_{\mathcal X}(t))-
 L_{\mathcal X}(t)>x}e^xdx,
}
where $ L_{\mathcal X}(t)$ is a pure jumps L\'{e}vy process
independent of $B$ with jump distribution $\mathcal X$.
The following lemma guaranties finiteness and positivity of this constant under
\textbf{A}.
\begin{lemma}\label{lemma_compound_poisson}
For any positive distribution $\mathcal X$ it holds that
$\mathcal H_{\mathcal X} \le 1$. Moreover, if $\mathcal X$ has finite variance,
then $\mathcal H_{\mathcal X} \in (0,1]$.
\end{lemma}

The next theorem gives the asymptotics of $\fr$ for all
$H\in (0,1)$.
\begin{theo} \label{infinite-horizon_theorem}
Assume that \textbf{A} holds if $H\le 1/2$. Then as $u \to \IF$
\bqn{\label{infinite-horizon_theorem_claim} \ \ \ \ \ \ \ \
\fr \sim
\begin{cases}
\mathbb{H}_{2H}\frac{2^{\frac{1}{2}-\frac{1}{2H}}\sqrt \pi}
{H^{1/2}(1-H)^{1/2}}
(C_Hu^{1-H})^{1/H-1}\Psi(C_Hu^{1-H}), & H>1/2 ,\\
\mathcal{H}_{2c^2\mathcal Z}  e^{-2cu}, & H=1/2, \\
\frac{\sqrt{2\pi} H^{H+1/2} u^H}{\mu c^{H+1}(1-H)^{H+1/2}}
\Psi(C_Hu^{1-H}), & H<1/2
\end{cases}
}
with $\mathcal{H}_{2c^2\mathcal Z} \in (0,1]$.
\end{theo}

\begin{remark}
Comparing the asymptotics of $\psi_{\delta,H}(u)$ in Proposition
\ref{proposition} and $\fr$
providing both asymptotics hold
we observe that as $u \to \IF$: \\
1) $\fr \sim \psi_{\delta,H}(u)$ for $H>1/2$ and all $\delta\ge 0$; \\
2) $\pi_{\mathcal Z,\frac{1}{2}}(u) \sim C_{\delta,\mathcal Z}
\psi_{\delta,\frac{1}{2}}(u)$ for all $\delta\ge 0$ with some
constant $C_{\delta,\mathcal Z}>0$ that
depends only on $\delta$ and $\mathcal Z$;\\
3.1) $\fr \sim \psi_{\mu,H}(u)$ for $H<1/2$;\\
3.2) $\fr = o(\psi_{H}(u))$ for $H<1/2$.
\end{remark}

A natural question is what happens if
condition \textbf{A} does not hold for $H\le 1/2$.
It seems to be very difficult to derive the exact asymptotics of
$\pi_{H,\mathcal Z}(u)$ in this case.
If  $H<1/2$ and $\mathbb{E}\{Z_1\}=\IF$
we give the following proposition:
\begin{sat}\label{proposition_infinite_expectation}
 Let $\mathbb{E}\{Z_1\}=\IF$ and $H<1/2$. Then as $u \to
\IF$
\bqn{\label{upper_bound_heavy_tail_H<1/2}
\fr = o(u^H \Psi(C_Hu^{1-H})),\\
\label{lower_bound_heavy_tail_H<1/2}
\Psi(C_Hu^{1-H}))= o(\fr).
}
\end{sat}


\section{Proofs}\label{sec_proofs}

\textbf{Proof of Theorem \ref{infinite-horizon_theorem}.}
\underline{Case $H>1/2$.} Upper bound.
Since all trajectories of $L(t), t\ge 0$ lay in $[0,\IF)$ we obtain
for all $u\ge 0$ $$\fr \le \psi_H(u).$$
Lower bound. As follows from
Corollary 3.2 in \cite{ParisRuinGenrealHashorvaJiDebicki}
\bqn{\label{proof_main_theo_paris}
\pk{\sup\limits_{t\ge 0}\inf\limits_{s\in [t,t+T_u]}(B(t)-ct)>u}
\sim \psi_H(u), \ \ \  u \to \IF}
if $T_u \le u^{2-1/H}a_u$ with any $a_u \to 0$ as $u \to \IF$.
Take in the following $T_u = u^{2-1/H}(\ln u)^{-1},u>1$.
We consider the subset $\mathcal A(u)$ of the general probability space
consisting of all outcomes such that
$\exists t\ge 0: B_H(t)-ct>u$, the subset
 $\mathcal A_{\inf}(u)$ consisting of
all outcomes such that $\exists t\ge 0: \inf\limits_{s\in [t,t+T_u]}(
B_H(s)-cs)>u$ and
 $\mathcal A_{L}(u)$ the subset of $\mathcal A_{\inf}(u)$ consisting
of all outcomes such that $\exists t\ge 0: B_H(L(t))-cL(t)>u$.
Denote by $G_u(t)$ the distribution function of
$$ \inf t\ge 0: \inf\limits_{s\in [t,t+T_u]}(B_H(s)-cs)>u
\Big|\exists z\ge 0: \inf\limits_{s_1\in [z,z+T_u]}(B_H(s_1)-cs_1)>u.$$
We have
\bqny{
\mathbb{P}(\mathcal A_{\inf}(u)\backslash \mathcal A_{L}(u))
&=&
\pk{\exists t\ge 0: \inf\limits_{s\in [t,t+T_u]}(B_H(s)-cs)>u,
\forall s_2 \ge0 : L(s_2) \notin [t,t+T_u]
}
\\&=&
\int\limits_0^\IF
\mathbb{P}\Big(\inf\limits_{s\in [t,t+T_u]}(B_H(s)-cs)>u,
\forall s_2 \ge0 : L(s_2) \notin [t,t+T_u]
\\& \ &
\Big|\exists z\ge 0: \inf\limits_{s_1\in [z,z+T_u]}(B_H(s_1)-cs_1)>u
\Big)
dG_u(t)
\\&=&
\int\limits_0^\IF
\mathbb{P}\Big(\inf\limits_{s\in [t,t+T_u]}(B_H(s)-cs)>u
\Big|\exists z\ge 0: \inf\limits_{s_1\in [z,z+T_u]}(B_H(s_1)-cs_1)>u
\Big)
\\& \ & \times \pk{\forall s_2 \ge0 : L(s_2) \notin [t,t+T_u]}
dG_u(t)
\\&\le&
\sup\limits_{t\ge 0} \pk{\forall s_2 \ge0 : L(s_2) \notin [t,t+T_u]}
\\& \ & \times
\int\limits_0^\IF
\mathbb{P}\Big(\inf\limits_{s\in [t,t+T_u]}(B_H(s)-cs)>u
\Big|\exists z\ge 0: \inf\limits_{s_1\in [z,z+T_u]}(B_H(s_1)-cs_1)>u
\Big)dG_u(t)
\\&=& \pk{\mathcal A_{\inf}(u)}
\sup\limits_{t\ge 0} \pk{\forall s \ge0 : L(s) \notin [t,t+T_u]}
.}

Finally,
\bqn{\label{generalization_Zu,H>1/2}
\fr\ge \pk{\mathcal A_{L}(u)} &=& \pk{\mathcal A_{\inf}(u)}
-\pk{\mathcal A_{\inf}(u)\backslash\mathcal A_{L}(u)}
\notag\\&\ge&
\pk{\mathcal A_{\inf}(u)}(1-\sup\limits_{t\ge 0}
\pk{\forall s \ge0 : L(s) \notin [t,t+T_u]})
\notag\\&=&
(1+o(1))\psi_H(u)(1-\sup\limits_{t\ge 0}
\pk{\forall s \ge0 : L(s) \notin [t,t+T_u]}),}
where the last line above follows by \eqref{proof_main_theo_paris}.
Since $Z_i$ are i.i.d. and $T_u \to \IF$ as $u\to \IF$,
the $\sup$ above tends to $0$ as $u$ tends to $\IF$ and the claim follows.
\\

\underline{Cases $H=1/2$ and $H<1/2$.}
We have by the self-similarity of fBm
with $t_0 = \frac{H}{c(1-H)}$ being the unique
maxima of $\Var(\frac{B_H(t)}{1+ct}), V(t) = \frac{B_H(t)}{1+ct},
\widetilde L_u(t) = L(t)/u$ and
$I_u(t_0) = [-u^{H-1}\ln u+t_0, t_0+u^{H-1}\ln u ]$
\bqn{\label{borel}
\fr & = &
\pk{\exists t\ge 0: B_H(u\widetilde L_u(t))>u+uc\widetilde L_u(t)}
\notag\\ &=& \pk{\exists t\ge 0:
\frac{B_H(\widetilde L_u(t))}{1+c\widetilde L_u(t)}>u^{1-H}}
\notag\\&\sim&
p(u) := \pk{\exists t \ge 0: \widetilde L_u(t)\in I_u(t_0)
, V(\widetilde L_u(t))> u^{1-H}}, \ \ \ \ u \to \IF.
}
The last line above follows from
Borell-TIS inequality (see \cite{20lectures}) and the asymptotics of
$p(u)$ given in \eqref{p(u)_h=1/2_asympt} for case $H=1/2$ and in
\eqref{H<1/2_finall_asympt} for case
$H<1/2$, see the Appendix for
the detailed proof. Now we treat cases $H=1/2$ and $H<1/2$ separately.
\newline

\underline{Case $H=1/2$.} We apply approach from \cite{HashorvaGrigori}.
Define $N_{u,S} = [\frac{\sqrt u \ln u}{S}]$ for $S>0$
and
$$c_{u,j,S} = t_0+jS/u, \ \ \ \Delta_{u,j,S} = [c_{u,j,S},c_{u,j+1,S}],$$
for $j\in [-N_{u,S},N_{u,S}]$.
We have for $S,u>0$ by Bonferroni inequality
\bqn{& \ &\notag
\sum\limits_{j=-N_{u,S}}^{N_{u,S}}
\pk{\exists t\ge0: \widetilde L_u(t) \in \Delta_{u,j,S} :
V(\widetilde L_u(t))>\sqrt u} - \Sigma(u,S)
\\&\le& p(u)\label{p(u)_h=1/2}
\\&\le&\notag
\sum\limits_{j=-N_{u,S}-1}^{N_{u,S}+1}
\pk{\exists t\ge0: \widetilde L_u(t) \in \Delta_{u,j,S},
V(\widetilde L_u(t))>\sqrt u}
,}
where
\bqny{
\Sigma(u,S) = \!\!\!\sum\limits_{-N_{u,S}-1\le i<j\le N_{u,S}+1}
\!\!\!\!\!\!\!\!\!
\pk{\exists t,s\ge0: \widetilde L_u(t) \in \Delta_{u,j,S},
\widetilde L_u(s) \in \Delta_{u,i,S},
V(\widetilde L_u(t))>\sqrt u, V(\widetilde L_u(s))>\sqrt u}.
}
Next we have as $u \to \IF$ and then $S\to \IF$
\bqn{\notag\label{double_sum}
\Sigma(u,S) &\le&\sum\limits_{-N_{u,S}-1\le i<j\le N_{u,S}+1}
\pk{\exists t,s \in \Delta_{u,j,S}\times \Delta_{u,i,S}:
V(t)>\sqrt u, V(s)>\sqrt u}
\\&=&o(e^{-2cu}),
}
where the last line above follows from \cite{Rolski17} eq. (43).
Next we approximate each probability in the sum above. We have
with $v=\sqrt u$ and $\phi_{u,j,S}$ the density of $B(c_{u,j,S})$,
$\widetilde L_{u,j,S}(t) = L(t)-uc_{u,j,S}$
\bqn{
& \ & \pk{\exists t \ge0: \widetilde L_u(t)
\in \Delta_{u,j,S}, \frac{B(\widetilde L_u(t))}
{1+c\widetilde L_u(t)}>\sqrt u}
\notag\\& = &
\pk{\exists t \ge0: \widetilde L_u(t)\in \Delta_{u,j,S},
B(\widetilde L_u(t))-vc\widetilde L_u(t)>v}
\notag\\& = &
\frac{1}{v}\int\limits_{\R}
\mathbb{P}\Big(
\exists t \ge 0: \widetilde L_u(t) \in \Delta_{u,j,S},
B(\widetilde L_u(t))-B(c_{u,j,S})
+B(c_{u,j,S})-vc(\widetilde L_u(t)-c_{u,j,S})-vcc_{u,j,S}>v
\notag\\& \ &
|B(c_{u,j,S}) = v-\frac{x}{v} \Big)
\phi_{u,j,S}(v-\frac{x}{v})dx
\notag\\&=&
\frac{1}{v}\int\limits_{\R}
\pk{\exists t\ge 0:  \widetilde L_u(t)\in \Delta_{u,j,S},
B(\widetilde L_u(t)-c_{u,j,S})+v-\frac{x}{v}
-vc(\widetilde L_u(t)-c_{u,j,S})-vcc_{u,j,S}>v}
\notag\\& \ & \times \phi_{u,j,S}(v-\frac{x}{v})dx
\notag\\&=&
\frac{1}{v}\int\limits_{\R}
\pk{\exists t \ge 0:L(t)\in u\Delta_{u,j,S},
 B(\widetilde L_{u,j,S}(t))-c\widetilde L_{u,j,S}(t)>
 x+v^2c c_{u,j,S}} \phi_{u,j,S}(v-x/v)dx
\notag\\&=&
\frac{1}{v}\int\limits_{\R} \pk{\exists
t\ge 0: \widetilde L_{u,j,S}(t) \in [0,S],
B(\widetilde L_{u,j,S}(t))-c\widetilde L_{u,j,S}(t)>
v^2 cc_{u,j,S}+x} \phi_{u,j,S}(v-x/v)dx
\notag\\& = &
\frac{1}{v}\int\limits_{\R}
\pk{\exists t\ge 0: \widetilde L_{u,j,S}(t) \in [0,S],
B(\widetilde L_{u,j,S}(t))-c\widetilde L_{u,j,S}(t)>x}
 \phi_{u,j,S}(v(1+cc_{u,j,S})-x/v)dx
\notag\\ & & \label{former_lemma}
\phantom{ghjl} \sim \frac{e^{-\frac{v^2(1+cc_{u,j,S})^2}{2c_{u,j,S}}}}
{\sqrt {2\pi c}}\frac{1}{v}
\int\limits_\R \pk{\exists t \ge0: L(t)\in [0,S],
B(L(t))-cL(t)>x}e^{2cx}dx,
}

as $u \to \IF$ and then $S \to \IF$,
we prove the last line above  in the Appendix.
The last integral above equals
\bqny{ & \ &
\frac{1}{2c}
\int\limits_{\R} \pk{\exists t\ge 0: 2c^2L(t) \in [0,2c^2S]
,\sqrt 2B(2c^2L(t))-2c^2L(t)>2cx}e^{2cx}d(2cx)
\\&=&
 \frac{cS}{2c^2S}
\int\limits_{\R} \pk{\exists t \ge 0: 2c^2L(t) \in [0,2c^2S],
\sqrt 2B(2c^2L(t))-2c^2L(t)>x}e^{x}dx
\\&\sim& cS\mathcal{H}_{2c^2L}, \ \ \ S \to \IF,
}
with $\mathcal{H}_{2c^2L}\in (0,1]$, where the last line above
follows from Lemma \ref{lemma_compound_poisson}.
 Summarizing all lines above we have
as $u \to \IF$ and then  $S \to \IF$
\bqny{
\sum\limits_{j=-N_{u,j,S}}^{N_{u,j,S}}
\pk{\exists t\ge0: \widetilde L_u(t) \in \Delta_{u,j,S} :
V(\widetilde L_u(t))>\sqrt u} &\sim&
\mathcal{H}_{2c^2 L}\sum\limits_{j=-N_u}^{j=N_u}
\frac{e^{-\frac{v^2(1+cc_j)^2}{2c_j}}}{\sqrt {2\pi c}v}Sc
\\&\sim& \mathcal{H}_{2c^2L}e^{-2cu}
}
where the last sum above was calculated in \cite{HashorvaGrigori}.
Combining the lines above with \eqref{p(u)_h=1/2} and \eqref{double_sum}
we have as $u\to \IF$
\bqn{\label{p(u)_h=1/2_asympt}
p(u) \sim \mathcal{H}_{2c^2L}e^{-2cu}}
and the claim follows.
\\

\underline{Case $H<1/2$.}
Let $\Omega_u = \#\{I_u(t_0)\cap \widetilde L_u(t)\}$. We show in the Appendix that
for sufficiently small $\varepsilon_0,\varepsilon_1,\varepsilon_2>0$
with $w = u^{\max(H+\varepsilon_1,1/H-2+\varepsilon_2)}$
as $u \to \IF$
\bqn{\label{H<1/2_capacity_of_main_interv}
p(u) \sim
\pk{\exists t\ge0:\widetilde L_u(t) \in I_u(t_0),V(\widetilde
 L_u(t))>u^{1-H}
\text{ and } \Omega_u \in [\varepsilon_0 u^H\ln u,w]}=:q(u).
}
We write first with $I_u^-(t_0) = [t_0-u^{H-1}\ln u,t_0]$ and
$I_u^+(t_0)=[t_0,t_0+u^{H-1}\ln u]$
\bqn{\label{positive_negative_intervals}
q(u) &\sim& \pk{\exists t \ge0: \widetilde L_u(t)\in I_u^-(t_0),
V(\widetilde L_u(t))>u^{1-H}\text{ and }\Omega_u \in
[\ve_0u^H\ln u,w]}
\\ & \ & +  \ \notag
\pk{\exists t\ge 0:\widetilde L_u(t)\in I_u^+(t_0),
V(\widetilde L_u(t))>u^{1-H} \text{ and } \Omega_u \in
[\ve_0u^H\ln u,w]}
\\ \notag &=:& q_-(u)+q_+(u),  \ \ u \to \IF,
}
proof is given in the Appendix.
Next we compute the asymptotic of $q_+(u)$.
Define for $i \ge 1$ $S_{u,i} = \frac{1}{u}\sum\limits_{j=1}
^{i}  Z_j$ ,
$t_u = \inf t\ge0:\widetilde\Pi_u(t)\ge t_0 $ and
$\xi_u = \widetilde\Pi_u(t_u)-t_0$.
We have by Bonferroni inequality
\bqn{& \ &\notag
\sum\limits_{i=1}^w \pk{V(t_0+S_{u,i})
>u^{1-H}\text{ and }\Omega_u \in [\ve_0 u^H\ln u,w]}
\\&\ge&\label{H<1/2_Bonferroni}
q_+(u)
\\&\ge&\notag
\sum\limits_{i=1 }^{\ve_0u^H\ln u} \pk{V(t_0+\xi_u+S_{u,i})
>u^{1-H}\text{ and }\Omega_u \in [\ve_0u^H\ln u,w]} -
\widetilde\Sigma(u),
}
where
$$\widetilde\Sigma(u) =
\sum\limits_{ 0\le j<i \le w} \pk{V(t_0+\xi_u+S_{u,i})
>u^{1-H}, V(t_0+\xi_u+ S_{u,j})>u^{1-H}} .$$

Thus, to compute the asymptotic of $q_+(u)$ we need to show, that
$\widetilde \Sigma(u)$
is negligible and calculate the asymptotics of both sums in
\eqref{H<1/2_Bonferroni}. We prove in the Appendix that
\bqn{\label{H<1/2_bound_Pi(u)}
\widetilde\Sigma(u) = o(u^H \Psi(C_Hu^{1-H})), \ \ \ u \to \IF.
}
\emph{Calculation of $\sum\limits_{i=1 }^{w} \pk{V(t_0+S_{u,i})
>u^{1-H}\text{ and }\Omega_u \in [\ve_0 u^H\ln u,w]}$.}
We can rewrite the sum in the upper bound in
\eqref{H<1/2_Bonferroni} as follows ($\mathcal N$ is a standard normal random
variable independent of all random variables
and stochastic processes which we consider):
\bqny{
\sum\limits_{i=1}^{w} \pk{\mathcal N>u^{1-H}
\frac{1+c(t_0+S_{u,i})}{(t_0+S_{u,i})^H}\text{ and }\Omega_u
\in [\ve_0u^H\ln u,w]}.}
 Setting
 $f(t) = \frac{1+ct}{t^H}$ we have $f(t_0)=C_H, f'(t_0) = 0$ and
$f''(t_0) = \frac{c^{H+2}(1-H)^{H+2}}{H^{H+1}}$.
Since for $t\in I(t_0)$ and any positive $\varepsilon$ for $u$ large enough
it holds that
$f(t_0)+\frac{(1-\varepsilon)f''(t_0)}{2}(t-t_0)^2<
f(t) < f(t_0)+\frac{(1+\varepsilon)f''(t_0)}{2}
(t-t_0)^2$ we have with $f^- = \frac{(1-\varepsilon)f''(t_0)}{2}$
\bqn{\label{H<1/2_quadratic_approximation}
\notag & \ &
\sum\limits_{0
\le i \le w} \pk{\mathcal N>
u^{1-H}C_H+u^{1-H}\frac{(1+\varepsilon)f''(t_0)}{2}S_{u,i}^2
\text{ and }\Omega_u \in [\ve_0u^H\ln u,w]}
\\&\le&
\sum\limits_{0
\le i \le w} \pk{\mathcal N>u^{1-H}
\frac{1+c(t_0+S_{u,i})}{(t_0+S_{u,i})^H}\text{ and }
\Omega_u \in [\ve_0u^H\ln u,w]}
\\&\le&\notag
\sum\limits_{0
\le i \le w} \pk{\mathcal N>u^{1-H}C_H+f^-u^{1-H}S_{u,i}^2
\text{ and }\Omega_u \in [\ve_0u^H\ln u,w]}
.}
Next we compute the asymptotic of the sum in the upper bound.
We have with $\widetilde N_{i}$ being an independent of $\mathcal N$
random variables with zero mean
$$S_{u,i} =
 \frac{(Z_1-\mu)+...+(Z_i-\mu)}{u}
+\frac{\mu i }{u} =: \frac{\widetilde N_{i}}{u}+\frac{\mu i}{u}
$$
and hence
\bqny{& \ &
\pk{\mathcal N>C_Hu^{1-H}+f^-u^{1-H}S_{u,i}^2
\text{ and }\Omega_u\in [\ve_0u^{H}\ln u,w]}
\\&=&
\pk{ \mathcal N>C_Hu^{1-H}+f^-u^{1-H}
(\frac{\widetilde N_{i}}{u}+\frac{\mu i}{u})^2
\text{ and }\Omega_u\in [\ve_0u^{H}\ln u,w]}.
}
Since for small $\varepsilon>0$ for all $i>u^\varepsilon$
very probably $\frac{\widetilde N_i}{u} << \frac{\mu i}{u}$ and very probably
$\Omega_u \in [\ve_0 u^H\ln u,w]$
we neglect term $\frac{\widetilde N_i}{u}$ and condition
$\Omega_u \in [\ve_0 u^H\ln u,w]$
in the sum above and write (proof is in the Appendix)
\bqn{\label{H<1/2_error_negligibility}\notag
& \ & \sum\limits_{0
\le i \le w } \pk{\mathcal N>
C_Hu^{1-H}+f^-u^{1-H}S_{u,i}^2\text{ and }\Omega_u\in
[\ve_0 u^{H}\ln u,w]}
\\&\sim&
\sum\limits_{0 \le i \le w}
\pk{ \mathcal N>C_Hu^{1-H}+f^-u^{1-H}(\frac{\mu i}{u})^2}
.}
For $i\ge \ve_0 u^H\ln u$ we have as $u \to \IF$
$$\pk{ \mathcal N>u^{1-H}f(t_0)+f^-u^{1-H}(\frac{\mu i}{u})^2}
\le \Psi(C_Hu^{1-H}+Cu^{-1-H}(u^H\ln u)^2) \sim
\Psi(C_Hu^{1-H})e^{-C\ln^2 u},
$$
hence
\bqn{\label{h<1/2_ln}
\sum\limits_{\ve_0u^H\ln u \le i \le w}
\pk{ \mathcal N>C_Hu^{1-H}+f^-u^{1-H}(\frac{\mu i}{u})^2} \le
w\Psi(C_Hu^{1-H})e^{-C\ln^2 u}
.}
Based on the proof of Theorem 2.2 in \cite{secondproj}
we have
\bqny{&\ &
\sum\limits_{0 \le i \le \ve_0 u^H\ln u}
\pk{ \mathcal N>C_Hu^{1-H}+f^-u^{1-H}(\frac{\mu i}{u})^2}
\\& = &
\sum\limits_{0 \le i \le \ve_0 u^H\ln u}
\Psi(C_Hu^{1-H}+f^-u^{1-H}(\frac{\mu i}{u})^2)
\\& = & (1+o(1))
\sum\limits_{0 \le i \le \ve_0 u^H\ln u}
\frac{(C_Hu^{1-H})^{-1}}{\sqrt{2\pi}}
e^{-(C_Hu^{1-H})^2/2-C_Hf^-u^{-2H}(\mu i)^2-
(f^-)^2u^{2-2H}(\frac{\mu i}{u})^4/2}
\\& = & (1+o(1))\Psi(C_Hu^{1-H})
\sum\limits_{0 \le i \le \ve_0 u^H\ln u}
  e^{- C_Hf^-u^{-2H}(\mu i)^2} ,
}

where the lines above follows
from the following inequality (see \cite{20lectures})
\bqn{\label{psi_asympt}
(1-\frac{1}{x^2})\frac{1}{\sqrt{2\pi}x}e^{-x^2/2}
\le \Psi(x)
\le \frac{1}{\sqrt{2\pi}x}e^{-x^2/2}, \ \ \ \ x>0,
}
 and the fact that
uniformly in $i\in [0,\ve_0 u^H\ln u]$ it holds that
$(f^-)^2u^{2-2H}(\frac{\mu i}{u})^4/2 \le u^{-1}$ as $u\to \IF$.
Next as $u\to \IF$
\bqny{ & \ &
\sum\limits_{0 \le i \le \ve_0u^H\ln u}
e^{- C_Hf^-u^{-2H}(\mu i)^2}
\\&=&
u^H\Big(\frac{1}{u^H}\sum\limits_{i \in [0, \ve_0 \ln u]\cap\Z u^{-H}}
e^{-i^2C_H\mu^2f^-}\Big)
\\&\sim&
u^H\int\limits_0^{\ve_0\ln u}
e^{-C_H\mu^2f^-z^2}dz
\\&\sim&\frac{u^H}{\sqrt{C_H\mu^2f^-}}
\int\limits_0^\IF
e^{-C_H\mu^2f^-z^2}d(z \sqrt{C_H\mu^2f^-})
\\&=&
\frac{\sqrt\pi u^H}{2\mu\sqrt{C_Hf^-}},
}
hence we obtain
$$
\sum\limits_{0 \le i \le \ve_0 u^{H}\ln u}
\pk{ \mathcal N>C_Hu^{1-H}+f^-u^{1-H}(\frac{\mu i}{u})^2}
\sim
\frac{\sqrt{2\pi}H^{H+1/2}}{2\mu\sqrt{1-\ve} c^{H+1}(1-H)^{H+1/2}}
u^H\Psi(C_Hu^{1-H}),\ u \to \IF.
$$
Letting
$\varepsilon \to 0$ and repeating all calculations for the lower bound
in \eqref{H<1/2_quadratic_approximation} by \eqref{h<1/2_ln} we obtain
\bqny{
\sum\limits_{i=1 }^{w} \pk{V(t_0+S_{u,i})
>u^{1-H}\text{ and }\Omega_u \in [\ve_0 u^H\ln u,w]} \sim
\frac{\sqrt{2\pi}H^{H+1/2}}{2\mu c^{H+1}(1-H)^{H+1/2}}
u^H\Psi(C_Hu^{1-H})
, \ \ \ u \to \IF.
}

\emph{Calculation of $\sum\limits_{i=1 }^{\ve_0 u^H\ln u}
\pk{V(t_0+\xi_u+S_{u,i})
>u^{1-H}\text{ and }\Omega_u \in [\ve_0 u^H\ln u,w]}$.} Since $\ve_0
u^H\ln u <w$, we compute
this sum by the same approach and finally we obtain that as $u \to \IF$
the sum is equivalent with
$$
\frac{\sqrt{2\pi}H^{H+1/2}}{2\mu c^{H+1}(1-H)^{H+1/2}}
u^H\Psi(C_Hu^{1-H}).
$$
Since the upper and the lower bounds in \eqref{H<1/2_Bonferroni}
are asymptotically equivalent we obtain by \eqref{H<1/2_bound_Pi(u)}
\bqny{
q_+(u)  \sim
\frac{\sqrt{2\pi}H^{H+1/2}}{2\mu c^{H+1}(1-H)^{H+1/2}}
u^H\Psi(C_Hu^{1-H}), \ \ \ u \to \IF.
}
Similarly we obtain that $q_-(u)$ has the same asymptotic,
and finally by the line above and
\eqref{positive_negative_intervals} we have
\bqn{\label{H<1/2_finall_asympt}
p(u) \sim \frac{\sqrt{2\pi}H^{H+1/2}}{\mu c^{H+1}(1-H)^{H+1/2}}
u^H\Psi(C_Hu^{1-H}), \ \ \ u \to \IF
}
and the claim follows by \eqref{borel}. \QED
\\

Before proving Lemma \ref{lemma_compound_poisson} we formulate
and prove the following auxiliary result.
\begin{lemma} \label{lem1}
For any $D(t),t\ge 0$
non-negative L\'{e}vy process independent of $B$ such that
$D(t) \to \IF$ as $\to \IF$ almost surely, it holds that
\bqn{\label{new_pick_const}
\limit{T} T^{-1}     \E{ \sup_{t\in [0,T]} e^{ B(D(t)) -
D(t)/2}} \in ( 0,\IF).
}
\end{lemma}

\textbf{Proof of Lemma \ref{lem1}.}
Let $Q(T,S)= \E{ \sup\limits_{t\in [T,S]} e^{ B(D(t)) - D(t)/2}}, S>T>0$.
Observe that for any non-negative random variable $K$ independent of $B$
 we have with $H_K$ distribution function of $K$
\bqn{\label{expectation}
\E{e^{B(K)-K/2}} = \int\limits_{\R}
\E{e^{B(x)-x/2}}dH_K(x) =
\int\limits_\R dH_K(x) = 1,
}
where we used that $\E{e^{B(x)-x/2}} = 1, x\ge 0$.
Next with $(B^*,D^*)$ independent copy of $(B,D)$
by the independence and stationarity of
the increments properties we have
\bqn{\label{ika}
Q(S,S+T)
&=& \E{\sup\limits_{t\in [S,S+T]}e^{B(D(t))-B(D(S))-\frac{D(t)-D(S)}{2}
+B(D(S))-D(S)/2}}
\notag\\&=&
\E{\sup\limits_{t\in [0,T]}e^{B(D(t))-D(t)/2}
e^{B^*(D^*(S))-D^*(S)/2}}
\notag\\&=&
\E{\sup\limits_{t\in [0,T]}e^{B(D(t))-D(t)/2}}
\E{e^{B^*(D^*(S))-D^*(S)/2}}
\notag\\&=&
Q(0,T),
}
where the last line above follows from \eqref{expectation}.
Let $Q(T) = Q(0,T)$. We have by the line above
\bqny{
Q(T+S)\le Q(0,T)+Q(T,T+S) = Q(T)+ Q(S).
}
Hence by Fekete's Lemma
$$ \limit{T} T^{-1} Q(T) = \inf_{T> 0} Q(T) /T \le Q(1) < \IF.$$
Our next aim is to show that the limit in \eqref{new_pick_const} is
positive.
We have by \eqref{ika}
\bqny{
\limit{T}  T^{-1}  Q(T) & =& \limit{T} T^{-1}
\sum_{j=0}^{T-1}(Q(T-j)- Q(T-j-1))
\\& = & \limit{T} T^{-1} \sum_{j=0}^{T-1} (Q(0, T-j)- Q(1,T-j))
\\&\ge & \limit{T} T^{-1}\sum_{j=0}^{T-1}
\E{  \max(0,  1- \sup_{t\in   [1,T-j] }
e^{ B(D(t)) - D(t)/2})}
\\&\ge&
\limit{T} T^{-1}\sum_{j=0}^{T-1}
\E{  \max(0,  1- \sup_{t\in   [1,\IF) }
e^{ B(D(t)) - D(t)/2})}
\\ &=& \E{  \max(0,  1- \sup_{t\in  [1,\IF) }
e^{ B(D(t)) - D(t)/2})}>0,
}
and the claim follows.\QED
\\

\textbf{Proof of Lemma \ref{lemma_compound_poisson}.}
\emph{Finiteness of $\mathcal H_{\mathcal X}$}.
We have
\bqny{\mathcal H_{\mathcal X} = \lim\limits_{S\to \IF}\frac{1}{S}
\E{\sup\limits_{t\ge0}e^{(B(L_{\mathcal X}(t))-
L_{\mathcal X}(t)/2)\mathbb I
(L_{\mathcal X}(t) \in [0,S])}}
 \le \lim\limits_{S\to\IF} \frac{1}{S}
 \E{\sup\limits_{z\in [0,S]}e^{B(z)-z/2}}
= \mathbb H_1=1, }
where we used that $\mathbb H_1=1,$ see, e.g., \cite{DeM15}.\\

\emph{Positivity of $\mathcal H_{\mathcal X}$.}
Let $T^* = \inf t\ge 0:  L_\mathcal X(t)>T$.
Then for $\theta = \frac{1}{2\mu}>0$
\bqny{& \ &
\E{\sup\limits_{t\ge0}e^{(B(L_{\mathcal X}(t))-
L_{\mathcal X}(t)/2)\mathbb I
(L_{\mathcal X}(t) \in [0,T])}}
\\&=&
\E{\sup\limits_{t\in[0,T^*)}e^{B(L_{\mathcal X}(t))-
L_{\mathcal X}(t)/2}}
\\&=&
\int\limits_{T^*\ge T\theta}\sup\limits_{t\in[0,T^*)}e^{B(
L_\mathcal X(t))-L_{\mathcal X}(t)/2}
d\mathbb P +
\int\limits_{T^*< T\theta}\sup\limits_{t\in[0,T^*)}e^{B(
L_{\mathcal X}(t))-L_{\mathcal X}(t)/2}
d\mathbb P \\&\ge&
\int\limits_{T^*\ge T\theta}\sup\limits_{t\in[0,T\theta]}e^{B(
L_{\mathcal X}(t))-L_{\mathcal X}(t)/2} d\mathbb P 
\\&=&
\E{(T^*\ge T\theta)\sup\limits_{t\in[0,T\theta]}e^{B(
L_{\mathcal X}(t))-L_{\mathcal X}(t)/2}}
\\&=&
\E{\mathbb I \sup\limits_{t\in[0,T\theta]}e^{B(L_{\mathcal X}(t)) -L_{\mathcal X}(t)/2}}
-\E{\mathbb  I(T^*< T\theta)\sup\limits_{t\in[0,T\theta]}e^{B(
L_{\mathcal X}(t))-L_{\mathcal X}(t)/2}}
.}

By Lemma \ref{lem1} we have that
$$\E{  \sup\limits_{t\in[0,T\theta]}e^{B(
L_{\mathcal X}(t))-L_{\mathcal X}(t)/2}} \ge C T,
\ \ \ \ T\to\IF$$
and hence to prove the claim it is enough to show that
\bqn{\label{small_mean}
\E{\mathbb  I(T^*< T\theta)\sup\limits_{t\in[0,T\theta]}e^{B(
L_{\mathcal X}(t))-L_{\mathcal X}(t)/2}}
= O(T), \ \ \ T \to \IF.
}

We have that the expression above equals
\bqny{ & \ &
\sum\limits_{k=0}^\IF
\E{\mathbb  I(T^*< T\theta)\mathbb I(L_{\mathcal X}
(T\theta)\in [kT,(k+1)T])
\sup\limits_{t\in[0,T\theta]}e^{B(L_{\mathcal X}(t))-
L_{\mathcal X}(t)/2}}
\\&\le&
\sum\limits_{k=0}^\IF
\E{\mathbb  I(T^*< T\theta)\mathbb I(L_{\mathcal X}
(T\theta)\in [kT,(k+1)T])
\sup\limits_{t\in[0,L_{\mathcal X}(T\theta)]}e^{B(t)-t/2}}
\\&\le&
\sum\limits_{k=0}^\IF
\E{\mathbb  I(T^*< T\theta)\mathbb I(L_{\mathcal X}
(T\theta)\in [kT,(k+1)T])
\sup\limits_{t\in[0,(k+1)T]}e^{B(t)-t/2}}
\\&\le&
CT\sum\limits_{k=0}^\IF
(k+1)\E{\mathbb  I(T^*< T\theta)\mathbb I(
L_{\mathcal X}(T\ve)\in [kT,(k+1)T])},
}

where the last line above follows from Lemma \ref{lem1}.
Thus, to prove the lemma it is enough to show that the sum in the
expression above tends to $0$ as $T\to \IF$. The sum above equals
\bqny{& \ &
\sum\limits_{k=0}^\IF
(k+1)\pk{T^*< T\theta, L_{\mathcal X}(T\theta)\in [kT,(k+1)T]}
\\&=&
\pk{T^*< T\theta}+
\sum\limits_{k=1}^\IF
\pk{T^*< T\theta, L_{\mathcal X}(T\theta)\ge kT}
\\&\le&
\pk{L_{\mathcal X}(T\theta)< T}+
\sum\limits_{k=1}^\IF
\pk{L_{\mathcal X}(T\theta)\ge kT}.
}
Since $\theta\mu = 1/2$
for the first term above we have by Chebyshev's inequality
\bqny{
\pk{L_{\mathcal X}(T\theta)< T} &=&
\pk{L_{\mathcal X}(T\theta)-
\E{L_{\mathcal X}(T\theta)}<T-\E{
L_{\mathcal X}(T\theta)}} \\&=&
\pk{L_{\mathcal X}(T\theta)-\mu T\theta<T(1-\mu\theta)}
\le \frac{\Var(L_{\mathcal X}(T\theta) - \mu T\theta)}{T^2/4}
= \frac{4\ve \Var(Z_1)}{T}.
}
Next by Chebyshev's inequality
\bqny{
\sum\limits_{k=1}^\IF
\pk{L_{\mathcal X}(T\theta)\ge kT} &\le& \sum\limits_{k=1}^\IF
\pk{L_{\mathcal X}(T\theta)-\mu \theta T\ge T(k-\mu\theta)}
\\&\le&
\sum\limits_{k=1}^\IF
\frac{\Var(L_{\mathcal X}(T\theta)-\mu \theta T)}{T^2(k-\mu\theta)^2}
= \frac{1}{T}\sum\limits_{k=1}^\IF\frac{\theta \Var(Z_1)}{(k-1/2)^2} = C/T
}
and hence
\bqny{
\pk{L_{\mathcal X}(T\theta)< T}+
\sum\limits_{k=1}^\IF \pk{L_{\mathcal X}(T\theta)\ge k T} 
\le C/T \to 0, \ \ \ T\to \IF
}
and the claim follows. \QED
\\

\textbf{Proof of Proposition \ref{proposition_infinite_expectation}.}
First we show \eqref{upper_bound_heavy_tail_H<1/2}.
Recall, that $q_+(u)$ and $q_{-}(u)$ defined in
\eqref{positive_negative_intervals}.
As in the proof of Theorem \ref{infinite-horizon_theorem}
by the proof of \eqref{borel} we have as $u \to \IF$
\bqny{
\fr \le q_-(u)+q_+(u)+\Psi(C_Hu^{1-H})e^{-C\ln^2u}.
}
Next
\bqny{
q_+(u) \le \sum\limits_{i=1}^w \pk{V(t_0+S_{u,i})
>u^{1-H}\text{ and }\Omega_u \in [\ve_0 u^H\ln u,w]}.
}
Now we fix large $A>0$ and define for
$i \ge 1$ $Y_{i,A} = Z_i\mathbb I(Z_i<A),
S_{u,i,A}' = \frac{1}{u}\sum\limits_{j=1}^{i} Y_{j,A}$.
Note that $Y_{i,A}$ has finite expectation and variance,
$Y_{i,A}\le Z_i$ and $Y_{i,A}$'s are i.i.d.
Since the variance of $V(t)$ has the unique global maxima at $t_0$,
decreases over $[t_0,\IF)$ and
$S_{u,i}\ge S_{u,i,A}'$
 we have
$$\sum\limits_{i=1}^w \pk{V(t_0+S_{u,i})
\!>\!u^{1-H}\text{ and }\Omega_u \in [\ve_0u^H\ln u,w]}
\!\le\!
\sum\limits_{i=1}^w \pk{V(t_0+S_{u,i,A}')
\!>\!u^{1-H}\text{ and }\Omega_u \in [\ve_0u^H\ln u,w]}.
$$
The asymptotic of the last sum above is already calculated
in the proof of Theorem \ref{infinite-horizon_theorem},
thus we have with
$\bar C = \frac{\sqrt{2\pi}H^{H+1/2}}{2c^{H+1}(1-H)^{H+1/2}}$
$$q_+(u) \le
(1+o(1))\frac{\bar C}{\E{Y_{1,A}}}u^H\Psi(C_Hu^{1-H}), \ \ u \to \IF.
$$
Estimating similarly $q_-(u)$ and summarizing all we have
$$\fr \le \frac{2\bar C+o(1)}{\E{Y_{1,A}}}u^H\Psi(C_Hu^{1-H}), \ \ \ u \to \IF.$$
Since $\E{Z_1} = \IF$ we have $\lim\limits_{A\to \IF}\E{Y_{1,A}} = \IF$ and
thus letting $A \to \IF$ we obtain the claim.
\\

Next we show \eqref{lower_bound_heavy_tail_H<1/2}. Fix some large
$M$, our aim is to show
\bqn{\label{lower_bound_heavy_tail_H<1/2_bound}
\fr \ge M\Psi(C_Hu^{1-H}), \ \ \ u \to \IF.
}
We choose large number $A_M$ such that
$\pk{ \# \{ [ut_0-A_M,ut_0+A_M]\cap \{L(t)_{t\ge 0} \}>2M\}}\ge 1-1/M$.
Define $R_{u,M} = \# \{ [ut_0-A_M,ut_0+A_M]\cap \{L(t)_{t\ge0}\} \}$.
We have for any small $\ve> 0$ by Bonferroni inequality
\bqny{
\fr &\ge& \pk{\exists t \ge 0: \widetilde \Pi_u(t)
\in [t_0-A_M/u,t_0+A_M/u],V(t)>u^{1-H}
\text{ and }
R_{u,M}\in [2M,u^\ve]}
\\&\ge& 2M\inf\limits_{t \in [t_0-A_M/u,t_0+A_M/u]}
\pk{V(t)>u^{1-H}\text{ and }R_{u,M}\in [2M,u^\ve]}
\\& \ &
- \sum\limits_{ -u^{\ve}\le j<i \le u^{\ve} }
\pk{V(t_0+S_{u,i})
>u^{1-H},
V(t_0+ S_{u,j})>u^{1-H}}
\\&\ge&
(1+o(1))2M\Psi(C_Hu^{1-H})\pk{R_{u,M}\in [2M,u^\ve]}
\\& \ & -Cu^{2\ve}\sup_{-u^{\ve}\le j<i \le u^{\ve}}
\pk{V(t_0+S_{u,i})>u^{1-H},V(t_0+ S_{u,j})>u^{1-H}}
\\&\ge&
1,5M\Psi(C_Hu^{1-H})-Cu^{2\ve}\sup_{-u^{\ve}\le j<i \le u^{\ve}}
\pk{V(t_0+S_{u,i})>u^{1-H},V(t_0+ S_{u,j})>u^{1-H}}
.}
Next repeating proof of \eqref{H<1/2_bound_Pi(u)} we have
$$\pk{V(t_0+S_{u,i})>u^{1-H},V(t_0+ S_{u,j})>u^{1-H}} \le
C\Psi(C_Hu^{1-H})u^{2-1/H}(\ln u)^{1/H}, \ \ \ u \to \IF
$$
and hence taking sufficiently small $\ve$ we have
\bqny{
\fr \ge  1,5M\Psi(C_Hu^{1-H})-C\Psi(C_Hu^{1-H})u^{2-1/H+2\ve}(\ln u)^{1/H}
\ge M\Psi(C_Hu^{1-H})
}
establishing \eqref{lower_bound_heavy_tail_H<1/2_bound} and thus
\eqref{lower_bound_heavy_tail_H<1/2} holds. \QED
\\

\section{Appendix}

\textbf{Proof of \eqref{borel}.}
We have for all $u>0$
\bqny{ & \ &
\pk{\exists t \ge 0: \widetilde L_u(t)\in I_u(t_0)
, V(\widetilde L_u(t))> u^{1-H}}
\\&\le&
\pk{\exists t\ge 0: V(\widetilde L_u(t))>u^{1-H}}
\\&\le&
\pk{\exists t \ge 0: \widetilde L_u(t)\in I_u(t_0),
V(\widetilde L_u(t))>u^{1-H}} +
\pk{\exists t \ge 0: \widetilde L_u(t)\notin I_u(t_0),
V(\widetilde L_u(t))>u^{1-H}}.
}
Next by Borell-TIS inequality (see \cite{20lectures}) we
have as $u \to \IF$
\bqny{\pk{\exists t \ge 0: \widetilde L_u(t)\notin I_u(t_0),
V(\widetilde L_u(t))>u^{1-H}}
\le \pk{\exists t\notin I_u(t_0): V(t)>u^{1-H}}
\le \Psi(C_Hu^{1-H})e^{-C\ln^2 u}.
}
Since
$\Psi(C_Hu^{1-H})e^{-C\ln^2 u} = o(p(u))$ as $u\to \IF$, see
\eqref{p(u)_h=1/2_asympt} for $H=1/2$ and \eqref{H<1/2_finall_asympt}
for $H<1/2$ the claim follows. \QED
\\

\textbf{Proof of \eqref{former_lemma}.}
We have
\bqny{& \ &
\int\limits_{\R}\pk{\exists t\ge 0: \widetilde L_{u,j,S}(t) \in [0,S],
B(\widetilde L_{u,j,S}(t))-c\widetilde L_{u,j,S}(t)>x}
\phi_{u,j,S}(v(1+cc_{u,j,S})-x/v)dx
\\&=&
\frac{e^{-\frac{v^2(1+cc_{u,j,S})^2}{2c_{u,j,S}}}}{\sqrt {2\pi c_{u,j,S}}}
\int\limits_{\R}\pk{\exists t\ge 0: \widetilde L_{u,j,S}(t) \in [0,S],
B(\widetilde L_{u,j,S}(t))-c\widetilde L_{u,j,S}(t)>x}
e^{\frac{x(1+cc_{u,j,S})}{c_{u,j,S}}-
x^2/(2uc_{u,j,S})} dx
}

and hence to prove the claim we need to show that as
$u \to \IF$ and then $S\to \IF$
\bqn{\label{app_h=1/2}& \ &
\int\limits_{\R}\pk{\exists t\ge 0: \widetilde L_{u,j,S}(t) \in [0,S],
B(\widetilde L_{u,j,S}(t))-c\widetilde L_{u,j,S}(t)>x}
e^{x(1+cc_{u,j,S})/c_{u,j,S}-
x^2/(2uc_{u,j,S})
} dx
\\&\sim&\notag
 \int\limits_\R \pk{\exists t \ge0: L(t)\in [0,S],
B(L(t))-cL(t)>x}e^{2cx}dx.
}

We have for all $j\in [-N_{u,S},N_{u,S}]$ with $M_u = u^{1/6}$
\bqny{& \ &
\int\limits_{\R}\pk{\exists t\ge 0: \widetilde L_{u,j,S}(t) \in [0,S],
B(\widetilde L_{u,j,S}(t))-c\widetilde L_{u,j,S}(t)>x}
e^{x(1+cc_{u,j,S})/c_{u,j,S}-x^2/(2uc_{u,j,S})} dx
\\&-& \int\limits_{\R}\pk{\exists t\ge 0: \widetilde
L_{u,j,S}(t) \in [0,S],B(\widetilde L_{u,j,S}(t))-
c\widetilde L_{u,j,S}(t)>x}e^{2cx} dx
\\&=&
\int\limits_{\R}\pk{\exists t\ge 0: \widetilde L_{u,j,S}(t) \in [0,S],
B(\widetilde L_{u,j,S}(t))-c\widetilde L_{u,j,S}(t)>x}
\left(e^{\frac{x(1+cc_{u,j,S})}{c_{u,j,S}}-x^2/(2uc_{u,j,S})}-
e^{2cx}\right) dx
\\&\le&
\int\limits_{\R}\pk{\exists t \in [0,S]:
B(t)-ct>x}e^{2cx}\left|e^{-cjSx/(uc_{u,j,S})-x^2/(2uc_{u,j,S})}-1\right| dx
\\&=& \int\limits_{-M_u}^{M_u}\pk{\exists t \in [0,S]:
B(t)-ct>x}e^{2cx}\left|e^{-cjSx/(uc_{u,j,S})-x^2/(2uc_{u,j,S})}-1\right| dx
\\ & \ &+\int\limits_{|x|>M_u}\pk{\exists t \in [0,S]:
B(t)-ct>x}e^{2cx}\left|e^{-cjSx/(uc_{u,j,S})-x^2/(2uc_{u,j,S})}-1\right| dx
\\&=:& I_1+I_2.
}

We have for $x\in [-M_u,M_u]$
\bqny{
|e^{-cjSx/(uc_{u,j,S})-x^2/(2uc_{u,j,S})}-1|
\le \sup\limits_{x\in [-M_u,M_u]}
|\frac{cjSx}{uc_{u,j,S}}+\frac{x^2}{2uc_{u,j,S}}|\le
\frac{C_1M_u\sqrt u \ln u+C_2M_u^2}{u}\le Cu^{-1/4}
}
and hence with some constants $C_S>0$ that does not depend on $u$
\bqn{\label{app_I1}
\notag I_1 &\le&
Cu^{-1/4}\int\limits_{-M_u}^{M_u}\pk{\exists t \in [0,S]:
B(t)-ct>x}e^{2cx}dx
\\&\le& \int\limits_\R\pk{\exists t \in [0,S]:
B(t)-ct>x}e^{2cx}dx
\le C_Su^{-1/4},
}
where finiteness of the integral above follows, e.g., from
Borell-TIS inequality. Next we have
\bqny{
I_2&\le& \int\limits_{x>M_u}\pk{\exists t \in [0,S]:
B(t)-ct>x}e^{2cx}\left(e^{cjSx/(uc_{u,j,S})-x^2/(2uc_{u,j,S})}-1\right) dx
+
\int\limits_{x<-M_u}e^{2cx} dx
\\&\le&
\int\limits_{x>M_u}\pk{\exists t \in [0,S]:
B(t)>x}e^{2cx}dx+ e^{-2cM_u}/2c
\\&\le&
\int\limits_{x>M_u}e^{-Cx^2/S+2cx}dx+ e^{-2cM_u}/2c,
}
where the line above follows from Borell-TIS inequality.
Combining the line above with \eqref{app_I1}
we obtain that as $u \to \IF$ and then $S\to \IF$
\bqny{& \ &
\int\limits_{\R}\pk{\exists t\ge 0: \widetilde L_{u,j,S}(t)\in [0,S],
B(\widetilde L_{u,j,S}(t))-c(\widetilde L_{u,j,S}(t))>x}
e^{x(1+cc_{u,j,S})-x^2/(2uc_{u,j,S})} dx
\\&\sim& \int\limits_{\R}
\pk{\exists t\ge 0: \widetilde L_{u,j,S}(t) \in [0,S],
B(\widetilde L_{u,j,S}(t))-c\widetilde L_{u,j,S}(t)>x}e^{2cx} dx.}

Denote $X_{u,j,S} = \inf t\ge 0: L(t)-uc_{u,j,S}\ge 0$ and
$\xi_{u,j,S} = L(X_{u,j,S}) - uc_{u,j,S}$. Note that $X_{u,j,S}$ is
a well defined random variable while $\xi_{u,j,S}$ is a non-negative
random variable.
Since the distribution of $L(t) - L(X_{u,j,S})$ for $t\ge X_{u,j,S}$
coincides with distribution of $L(t)$ for $t\ge 0$ we have
\bqny{
& \ &
\pk{\exists t\ge 0: \widetilde L_{u,j,S}(t) \in [0,S],
B(\widetilde L_{u,j,S}(t))-c\widetilde L_{u,j,S}(t)>x}
\\&=&
\pk{\exists t\ge X_{u,j,S}: \widetilde L_{u,j,S}(t) \in [0,S],
B(\widetilde L_{u,j,S}(t))-c\widetilde L_{u,j,S}(t)>x}
\\&=&
\mathbb P \{\exists t\ge X_{u,j,S}: (L(t)-L(X_{u,j,S}))
+(L(X_{u,j,S})-uc_{u,j,S}) \in [0,S],
\\& \ &
B\Big(
(L(t)-L(X_{u,j,S}))
+(L(X_{u,j,S})-uc_{u,j,S})
\Big)
\\& \ &-c\Big((L(t)-L(X_{u,j,S}))
+(L(X_{u,j,S})-uc_{u,j,S})\Big)>x\}
\\&=&
\pk{\exists z\ge0: \hat L(z)
+\xi_{u,j,S} \in [0,S],
B(\hat L(z)+\xi_{u,j,S})-c(\hat L(z)+\xi_{u,j,S})>x},
}
where $\hat L$ is an independent of $\xi_{u,j,S}$ copy of $L$.
We have with $F_{\xi_{u,j,S}}$ the df of $\xi_{u,j,S}$,
$\phi_a$ density of $B(a)$
that the last probability above equals

\bqny{
& \ &
\int\limits_{0}^S \pk{\exists t\ge0: L(t)
+a \in [0,S], B(L(t)+a)-c(L(t)+a)>x} dF_{\xi_{u,j,S}}(a)
\\& = &
\int\limits_{0}^S \pk{\exists t\ge0: K(t) \in [0,S-a],
B(L(t)+a)-B(a)-cL(t)>x+ac-B(a)} dF_{\xi_{u,j,S}}(a)
\\& = &
\int\limits_{0}^S \int\limits_\R\pk{\exists t\ge0: L(t) \in [0,S-a],
B(L(t))-cL(t)>x+ac-y} \phi_a(y)dydF_{\xi_{u,j,S}}(a)
\\& = &
\int\limits_{0}^S\int\limits_\R \pk{\exists t\ge0: L(t) \in [0,S-a],
B(L(t))-cL(t)>z} \phi_a(x+ac-z)dzdF_{\xi_{u,j,S}}(a)
\\& =: &
\frac{1}{\sqrt{2\pi}}\int\limits_{0}^S \int\limits_\R
W(a,z)e^{-\frac{(x+ac-z)^2}{2a}}a^{-1/2} dzdF_{\xi_{u,j,S}}(a)
}

and hence
\bqny{& \ &
\int\limits_{\R}\pk{\exists t\ge 0: \widetilde L(t) \in [0,S],
B(\widetilde L(t))-c\widetilde L(t)>x}e^{2cx} dx
\\&=&
\frac{1}{\sqrt{2\pi}}\int\limits_\R\int\limits_{0}^S\int\limits_\R
W(a,z)e^{-\frac{(x+ac-z)^2}{2a}+2cx}a^{-1/2}
dzdF_{\xi_{u,j,S}}(a)dx.
}
Changing variable $x\to w+z-ac$ we have that the last expression above equals

\bqny{& \ &
\frac{1}{\sqrt{2\pi}}\int\limits_\R\int\limits_{0}^S\int\limits_\R
W(a,z)e^{-\frac{w^2}{2a}+2c(w-ac)}a^{-1/2}e^{2cz}
dzdF_{\xi_{u,j,S}}(a)dw
\\&=&
\int\limits_\R\int\limits_{0}^S
\frac{e^{-\frac{(w-2ac)^2}{2a}}}{\sqrt{2\pi a}}
\int\limits_\R
W(a,z)e^{2cz}
dzdF_{\xi_{u,j,S}}(a)dw
\\&=&
\int\limits_{0}^S\int\limits_\R
W(a,z)e^{2cz}
dzdF_{\xi_{u,j,S}}(a),
}
where we used that
$\int\limits_\R e^{-(w-2ac)^2/2a}dw/\sqrt{2\pi a}=1$ for
$a>0$. Next we give the bounds for the integral above. We have
for the upper bound
\bqny{& \ &
\int\limits_{0}^S\int\limits_\R
W(a,z)e^{2cz}
dzdF_{\xi_{u,j,S}}(a)
\\&\le&
\int\limits_{0}^S\int\limits_\R
W(0,z)e^{2cz}
dzdF_{\xi_{u,j,S}}(a)
\\&=&
\pk{\xi_{u,j,S}\le S} \int\limits_\R
\pk{\exists t\ge0: L(t) \in [0,S],
B(L(t))-cL(t)>z}e^{2cz}dz.
}
Lower bound. Fix any small $\ve>0$.
\bqny{
& \ &
\int\limits_{0}^S\int\limits_\R
W(a,z)e^{2cz}dzdF_{\xi_{u,j,S}}(a)
\\ &\ge &
\int\limits_0^{\ve S}\int\limits_\R
W(a,z)e^{2cz}dzdF_{\xi_{u,j,S}}(a)
\\ &\ge &
\pk{\xi_{u,j,S}\le\ve S}\int\limits_\R
\pk{\exists t\ge0: L(t) \in [0,S(1-\ve)],
B(L(t))-cL(t)>z}e^{2cz}dz.
}
We have
$\pk{\xi_{u,j,S}>\ve S} \le \pk{Z_1>\ve S}$ and thus uniformly for $u$
\bqn{\label{E_{u,S}_asympt}
\pk{\xi_{u,j,S}>\ve S} \to 0, \ \ \ S \to \IF.
}
Thus, as $u\to \IF$ and then $S\to \IF$ we have
\bqny{
& \ & (1+o(1)) \int\limits_\R
\pk{\exists t\ge0: L(t) \in [0,S(1-\ve)],
B(L(t))-cL(t)>z}e^{2cz}dz
\\&\le&
\int\limits_{\R}\pk{\exists t\ge 0: \widetilde L_{u,j,S}(t) \in [0,S],
B(\widetilde L_{u,j,S}(t))-c\widetilde L_{u,j,S}(t)>x}e^{2cx} dx
\\&\le&
(1+o(1)) \int\limits_\R
\pk{\exists t\ge0: L(t) \in [0,S],
B(L(t))-cL(t)>z}e^{2cz}dz.
}
Repeating the calculations for the lower bound above as in Theorem
\ref{infinite-horizon_theorem} and then letting $\ve \to 0$
we obtain that the bounds above asymptotically agree and the proof
is completed.
\QED

Before giving the proofs of the remaining parts we prove the following
lemma.
\begin{lemma} \label{lemma_sum_iid} Let $V_1,...,V_n$ be i.i.d. random
variables having finite variance and expectation $\mu_V$. Then for
any $\ve>0$ and $n\ge 1$
\bqny{
\pk{V_1+...+V_n>(\mu_V+\ve)n}\le \hat Cn^{-1},
\quad
\pk{V_1+...+V_n>(\mu_V-\ve)n}\le \hat C n^{-1}
}
with $\hat C = \Var(V_1)/\ve^2$.
\end{lemma}

\textbf{Proof of Lemma \ref{lemma_sum_iid}.}
We have
\bqny{
\pk{V_1+...+V_n>(\mu_V+\ve)n} &=& \pk{(V_1-\mu_V)+...+(V_n-\mu_V)
>n\ve} \\&\le& \frac{\Var((V_1-\mu_V)+...+(V_n-\mu_V))}{n^2\ve^2}
=\hat C n^{-1},
}
where the inequality above follows from Chebyshev's inequality.
The second claim of the lemma follows exactly by the same arguments.\QED
\\

\textbf{Proof of \eqref{H<1/2_capacity_of_main_interv}.}
Recall, $\Omega_u = \#\{I(t_0)\cap \widetilde L_u(t)\}$.
We have
\bqny{
\pk{\Omega_u>w} &\le& \pk{S_{u,i}<2u^{H-1}\ln u, i>w }
\\&\le& \pk{Z_1+Z_2+...+Z_{[w]}<2u^H\ln u }
\le
\pk{Z_1+Z_2+...+Z_{[w]}<\mu w/2 }
 \le Cw^{-1},
}
where the last line above follows from Lemma \ref{lemma_sum_iid}.
Next we have by the asymptotics of $\psi_H(u)$ given in Proposition
\eqref{proposition}
\bqn{\notag& \ &
\pk{\exists t \ge 0: \widetilde L_u(t)\in I_u(t_0),
V(\widetilde L_u(t))>u^{1-H}
\text{ and }\Omega_u>w}
\\&\le&\notag
\pk{\exists t \in I_u(t_0): V(t)>u^{1-H} \text{ and }
\Omega_u>w}
\\&\le&\notag
\psi_H(u)\pk{\Omega_u>w}
\\&\le&\notag
Cu^{(1/H-1)(1-H)}\Psi(C_Hu^{1-H})w^{-1}
\\&=&\notag
u^{1/H-2-\max(H+\varepsilon_1,1/H-2+\varepsilon_2)}
 C u^H\Psi(C_Hu^{1-H})
\\&=&\label{App_capacity_mainint_upperbound}
o(u^H\Psi(C_Hu^{1-H})), \ \ \ u \to \IF
.}
Since $\Omega_u$ is independent of $B_H$ we have for large $u$
by Lemma \ref{lemma_sum_iid}
\bqny{ & \ &
\pk{\exists t\ge0 :\widetilde
L_u(t) \in I_u(t_0), V(t)\!>\!u^{1-H}\text{ and }
\Omega_u<\ve_0 u^H\ln u}
\\&\le&
\pk{V(t_0)\!>\!u^{1-H}} Cu^H\ln u \pk{\Omega_u<\ve_0 u^H\ln u}
\\ &\le&
Cu^H\ln u\Psi(C_Hu^{1-H}) u^{-H}(\ln u)^{-1}
\\&=&C\Psi(C_Hu^{1-H}).
}
Thus the claim follows by the line above,
\eqref{App_capacity_mainint_upperbound} and
the asymptotic of $p(u)$ given in
\eqref{H<1/2_finall_asympt}.
\QED \\

\textbf{Proof of \eqref{H<1/2_bound_Pi(u)}.}
By Lemma 2.3 in \cite{PicandsA} we have for
$\mathcal N_1$ and $\mathcal N_2$ being standard Gaussian random variables
with correlation $r$
\bqn{\label{pick_inequality}
\pk{\mathcal N_1>x,\mathcal N_2>x} \le 3\Psi(x)\Psi(x\sqrt {\frac{1-r}{1+r}}),
 \ \ \ x>0.}
We analyze each summand of $\Sigma(u)$ by the inequality above.
Assume that $i-j = n_{i,j}>0$. We have with
$K_{u,j} = t_0+\xi_u+S_{u,j}$, $Z_{1}+...+Z_{n_{i,j}} = Z_{i,j}$,
$r_k(z) = \cor(V(k),V(k+z/u)), \ \sigma(t) = \Var(V(t)),
\bar V(t) = V(t)/\sigma(t), F_{K_{u,j}}$ the distribution function of
$K_{u,j}$ and $F_{Z_{i,j}}$ the distribution function of $Z_{i,j}$
\bqny{ & \ &
\pk{V(t_0+\xi_u+S_{u,j})
>u^{1-H},V(t_0+\xi_u+S_{u,i})>u^{1-H}}
\\&=&
\pk{V(K_{u,j})>u^{1-H}, V(K_{u,j}+Z_{i,j}/u)>u^{1-H}}
\\&=&
\int\limits_0^{\IF} \pk{\bar V(K_{u,j})>\frac{u^{1-H}}{\sigma(K_{u,j})},
\bar V(K_{u,j}+z/u) > \frac{u^{1-H}}{\sigma(K_{u,j}+\frac{z}{u})}}dF
_{Z_{i,j}}(z)
\\&\le&
\int\limits_0^{\IF} \pk{\bar V(K_{u,j})>\frac{u^{1-H}}{\sigma(t_0)},
\bar V(K_{u,j}+z/u)>\frac{u^{1-H}}{\sigma(t_0)}}dF_{Z_{i,j}}(z)
\\&=&
\int\limits_0^{\IF} \int\limits_{t_0}^{\IF}
\pk{\bar V(k)>\frac{u^{1-H}}{\sigma(t_0)},
\bar V(k+z/u)>\frac{u^{1-H}}{\sigma(t_0)}}dF_{K_{u,j}}(k)dF_{Z_{i,j}}(z)
\\&\le&
3\Psi(C_Hu^{1-H})\int\limits_{t_0}^{\IF}\int\limits_0^{\IF}
 \Psi(C_Hu^{1-H}\sqrt{\frac{1-r_k(z)}{2}})dF_{Z_{i,j}}(z)dF_{K_{u,j}}(k),
}
where we used that $K_{u,j}\ge t_0$ always and \eqref{pick_inequality}.
We have
$r_k(z) = corr(B_H(k),B_H(k+z/u)) = 1-\frac{1}{2t_0^{2H}}|z/u|^{2H}+
o(u^{-2H})$ uniformly for $z \in [0,u^{H}\ln u]$ and
$k\ge t_0>0$. Since $r_k(z)$ is decreasing for
$z\in [0,\IF)$ and fixed $k$ we have that the internal integral in the
line above does not exceed
\bqny{ & \ &
 \int\limits_0^{u^H\ln u }
 \Psi(Cu^{1-2H}z^H)dF_{Z_{i,j}}(z)
 + \int\limits_{u^H\ln u }^\IF
 \Psi(u^{1-H}\sqrt{\frac{1-r_k(z)}{2}}) dF_{Z_{i,j}}(z)
 \\&\le&
\int\limits_0^\IF
 \Psi(Cu^{1-2H}z^H)dF_{Z_{i,j}}(z)
 +  \Psi(u^{1-H}\sqrt{\frac{1-r_k(u^H\ln u)}{2}})
 \\  &\le&
 \int\limits_0^\IF
 \Psi(Cu^{1-2H}z^H)dF_{Z_{i,j}}(z)
 +  \Psi(u^{(1-H)^2/2})
}
and hence we have
\bqny{
\int\limits_{t_0}^{\IF}\int\limits_0^{\IF}
 \Psi(C_Hu^{1-H}\sqrt{\frac{1-r_k(z)}{2}})dF_{Z_{i,j}}(z)dF_{K_{u,j}}(k) \le
 \int\limits_0^\IF
 \Psi(Cu^{1-2H}z^H)dF_{Z_{i,j}}(z)
 +  \Psi(u^{(1-H)^2/2}) .
}

We have
\bqny{\notag
\int\limits_0^\IF \Psi(Cu^{1-2H}z^H) dF_{Z_{i,j}}(z)
&\le& \int\limits_0^{u^{2-1/H}(\ln u)^{1/H}}
 \Psi(Cu^{1-2H}z^H)dF_{Z_{i,j}}(z)+
 \int\limits_{u^{2-1/H}(\ln u)^{1/H}}^\IF
 \Psi(C\ln u)dF_{Z_{i,j}}(z)
 \\&\le&\notag
 C u^{2-1/H}(\ln u)^{1/H} F_{Z_{i,j}}(u^{2-1/H}(\ln u)^{1/H})+
 e^{-C\ln^2 u} 
\\&\le& C u^{2-1/H}(\ln u)^{1/H} \pk{Z_1+...+Z_{n_{i,j}}<u^{2-1/H}(\ln u)^{1/H}}
 + e^{-C\ln^2 u}
 \notag \\&\le& C u^{2-1/H}(\ln u)^{1/H}/(i-j)
  + e^{-C\ln^2 u},
}
where the last line above follows by Lemma \ref{lemma_sum_iid}.
Thus, for all $i,j \in [0,u^{\max(H+\varepsilon_1,1/H-2+\varepsilon_2)}]$
\bqny{
 \pk{V(t_0+\xi_u+S_{u,j})>u^{1-H},V(t_0+\xi_u+S_{u,i})>u^{1-H}}
\le
C\Psi(C_Hu^{1-H})u^{2-1/H}(\ln u)^{\frac{1}{H}}
(i-j)^{-1}.
}
Using the inequality above to prove the claim we need to show that
\bqn{\label{app_H<1/2_aux}
u^{2-1/H}(\ln u)^{1/H} \sum\limits_{0\le j<i \le
u^{\max(H+\varepsilon_1,1/H-2+\varepsilon_2)}}
(i-j)^{-1} = o(u^H), \ \ u \to \IF.
}
By the fact that $\sum\limits_{0\le j<i \le
m} (i-j)^{-1} \le C m\ln m$ for large $m$ we complete
the proof since for sufficiently
small $\ve_1$ and $\ve_2$ the expression above does not exceed
\bqny{
u^{2-1/H}(\ln u)^{1/H}
u^{\max(H+\varepsilon_1,1/H-2+\varepsilon_2)} C\ln u
= o(u^H), \ \ u \to \IF. \ \ \  \ \ \ \ \ \ \ \ \ \ \ \ \  \ \ \ \ \ \ \ \
 \ \ \
\hfill \Box}

\textbf{Proof of \eqref{positive_negative_intervals}.}
To prove the claim in view of the final asymptotic of
$p(u)$ given in \eqref{H<1/2_finall_asympt} we need to show that
as $ u \to \IF$
\bqny{ & \ &
\mathbb P \Big(\exists t_1,t_2\ge 0:\widetilde L_u(t_1) \in I_u^-(t_0),
\widetilde L_u(t_2) \in I_u^+(t_0),
V(\widetilde L_u(t_1))>u^{1-H},
V(\widetilde L_u(t_2))>u^{1-H}
\\& \ & \text{ \ \ \ \ \ and } \Omega_u \in [\ve_0 u^H\ln u,w]\Big)
\!= o(u^H\Psi(C_Hu^{1-H})).
}
Define for $j<0$ $S_{u,j} = -\frac{1}{u}\sum\limits_{k=1}^{|j|} Z'_k,$
where $Z_k'$ is an independent copy of $Z_k$. We have that the probability
above does not exceed
\bqny{
\sum\limits_{-w\le j\le 0 ,0 \le i\le w }
\pk{V(t_0+S_{u,j})
>u^{1-H},V(t_0+S_{u,i})>u^{1-H}}.
}
Repeating the proof of \eqref{H<1/2_bound_Pi(u)} we obtain the claim,
the only difference is that in \eqref{app_H<1/2_aux} the index of summation
changes to $-u^{\max(H+\varepsilon_1,1/H-2+\varepsilon_2)}\le j\le 0,0\le
i \le u^{\max(H+\varepsilon_1,1/H-2+\varepsilon_2)}$. \QED
\\

\textbf{Proof of \eqref{H<1/2_error_negligibility}.}
We have with $F_{\widetilde N_i}$ the df of $\widetilde N_{i}$
for some sufficiently small $\varepsilon>0$
\bqn{\label{ints}
& \ &
 \pk{\mathcal N> C_Hu^{1-H}+u^{1-H}f_-S_{u,i}^2\text{ and }
 \Omega_u\in [\ve_0 u^H\ln u,w]}
\notag\\&=&
 \pk{\mathcal N>C_Hu^{1-H}+u^{1-H}f_-(
\frac{\widetilde N_i}{u}+\frac{\mu i}{u}
)^2\text{ and }\Omega_u\in [\ve_0 u^H\ln u,w]}
\notag\\&=&
\int\limits_\R \pk{\mathcal N>C_Hu^{1-H}+u^{1-H}f_-(
\frac{z}{u}+\frac{\mu i}{u}
)^2\text{ and }\Omega_u\in [\ve_0 u^H\ln u,w]}dF_{\widetilde N_i}(z)
\notag\\&=&
\int\limits_{|z|\le i/\ln i} \pk{\mathcal N>C_Hu^{1-H}+u^{1-H}f_-(
\frac{z}{u}+\frac{\mu i}{u}
)^2\text{ and }\Omega_u\in [\ve_0 u^H\ln u,w]}dF_{\widetilde N_i}(z)
\notag\\& \ & +
\int\limits_{|z|> i/\ln i} \pk{\mathcal N>C_Hu^{1-H}+u^{1-H}
f_-(\frac{z}{u}+\frac{\mu i}{u}
)^2\text{ and }\Omega_u\in [\ve_0 u^H\ln u,w]}dF_{\widetilde N_i}(z).
}
The second integral does not exceed
\bqny{
\pk{\mathcal N>C_Hu^{1-H}}\pk{|\widetilde N_{i}|>i/\ln i}
\le C\Psi(C_Hu^{1-H})i^{-1}\ln^2 i,
}
where the last inequality follows from Chebyshev's inequality.
Thus,
\bqn{\notag & \ &
\sum\limits_{1  \le i \le w}
\int\limits_{|z|> i/\ln i} \pk{\mathcal N>C_Hu^{1-H}+u^{1-H}
f_-(\frac{z}{u}+\frac{\mu i}{u}
)^2\text{ and }\Omega_u\in [\ve_0 u^H\ln u,w]}dF_{\widetilde N_i}(z)
\\&\le&\notag
C\sum\limits_{1 \le i \le w}
\Psi(C_Hu^{1-H})i^{-1}\ln^2 i
\\&\le& \label{App_negligibility_big_z} C\ln^3 u\Psi(C_Hu^{1-H}).
}

For the first integral in \eqref{ints} we have for any small
$\bar\varepsilon_1>0$
\bqny{ & \ &
\int\limits_{|z|\le i/\ln i}
 \pk{\mathcal N>C_Hu^{1-H}+u^{1-H}f_-
( (1+\bar\varepsilon_1)\frac{\mu i}{u})^2\text{ and }
\Omega_u\in [\ve_0 u^H\ln u,w]}dF_{\widetilde N_i}(z)
\\&\le&
\int\limits_{|z|\le i/\ln i} \pk{\mathcal N>C_Hu^{1-H}+u^{1-H}
f_-(\frac{z}{u} +\frac{\mu i}{u})^2\text{ and }\Omega_u
\in [\ve_0 u^H\ln u,w]
}dF_{\widetilde N_i}(z)
\\ &\le&
\int\limits_{|z|\le i/\ln i}
 \pk{\mathcal N>C_Hu^{1-H}+u^{1-H}f_-
( (1-\bar\varepsilon_1)\frac{\mu i}{u})^2\text{ and }
\Omega_u\in [\ve_0 u^H\ln u,w]
}dF_{\widetilde N_i}(z).
}
The upper bound in the inequality above equals
\bqny{
 \pk{\mathcal N>C_Hu^{1-H}+u^{1-H}f_-
( (1-\bar\varepsilon_1)\frac{\mu i}{u})^2}(1-\pk{|\widetilde
N_{i}|>i/\ln i}) \pk{\Omega_u\in [\ve_0 u^H\ln u,w]}
.}

By Chebyshev's inequality and the proof of
\eqref{H<1/2_capacity_of_main_interv} we have
uniformly for $i>u^{H/2}$
$$1-u^{-H/4}\le 1-\pk{|\widetilde N_{i}|>i/\ln i}\le 1, \ \ \
1-u^{-H}\le \pk{\Omega_u\in [\ve_0u^H\ln u,w]}\le 1.$$
Thus, by the line above
\bqny{& \ &
\sum\limits_{1 \le i \le w}
\int\limits_{|z|\le i/\ln i}
 \pk{\mathcal N>C_Hu^{1-H}+u^{1-H}f_-(
\frac{z}{u}+
 \frac{\mu i}{u})^2\text{ and }\Omega_u\in
 [\ve_0 u^H\ln u,w]}dF_{\widetilde N_i}(z)
\\ &\le&
(1+o(1))\sum\limits_{u^{H/4} \le i \le w}
\pk{\mathcal N>C_Hu^{1-H}+u^{1-H}f_-
( (1-\bar\varepsilon_1)\frac{\mu i}{u})^2}+O(u^{H/4}\Psi(C_Hu^{1-H}))
\\&\le&
(1+o(1))\sum\limits_{0 \le i \le w}
\pk{\mathcal N>C_Hu^{1-H}+u^{1-H}f_-
( (1-\bar\varepsilon_1)\frac{\mu i}{u})^2}+Cu^{H/4}\Psi(C_Hu^{1-H})
.}

Similarly for the lower bound we have
\bqny{& \ &
\sum\limits_{1 \le i \le w}
\int\limits_{|z|\le i/\ln i}
 \pk{\mathcal N>C_Hu^{1-H}+u^{1-H}f_-(
\frac{z}{u}+
 \frac{\mu i}{u})^2\text{ and }\Omega_u\in [\ve_0 u^H\ln u,w]}
 dF_{\widetilde N_i}(z)
\\ &\ge&
(1+o(1))\sum\limits_{u^{H/4} \le i \le w}
\pk{\mathcal N>C_Hu^{1-H}+u^{1-H}f_-
( (1+\bar\varepsilon_1)\frac{\mu i}{u})^2}
\\ &\ge&
(1+o(1))\sum\limits_{0 \le i \le w}
\pk{\mathcal N>C_Hu^{1-H}+u^{1-H}f_-
( (1+\bar\varepsilon_1)\frac{\mu i}{u})^2}-Cu^{H/4}\Psi(C_Hu^{1-H})
.}

Finally by \eqref{App_negligibility_big_z} we have
\bqny{ & \ & (1+o(1))
\sum\limits_{0 \le i \le w}
\pk{\mathcal N>C_Hu^{1-H}+u^{1-H}f_-
( (1+\bar\varepsilon_1)\frac{\mu i}{u})^2}-Cu^{H/4}\Psi(C_Hu^{1-H})
\\&\le&
\pk{\mathcal N> C_H u^{1-H}+u^{1-H}f_-S_{u,i}^2\text{ and }
\Omega_u\in [\ve_0u^H\ln u,w]}
\\&\le&
(1+o(1))\!\!\sum\limits_{0 \le i \le w}\!\!\!
\pk{\mathcal N>C_Hu^{1-H}+u^{1-H}f_-
( (1-\bar\varepsilon_1)\frac{\mu i}{u})^2}+Cu^{\frac{H}{4}}\Psi(C_Hu^{1-H})+
C\ln^3 u\Psi(C_Hu^{1-H})
.}

Repeating the calculations for the sums above as in the proof of Theorem
\ref{infinite-horizon_theorem} after
\eqref{H<1/2_error_negligibility} and then letting $\varepsilon_1\to 0$ we
obtain the claim, since
$u^{H/4}\Psi(C_Hu^{1-H}) = o(p(u)), \ u \to \IF$
by \eqref{H<1/2_finall_asympt}.
\QED\\

\textbf{Acknowledgement.} The Author would like to thank Krzysztof
{D\c{e}bicki and Enkelejd Hashorva for numerous discussions on the
topic of the paper.
The Author was supported by the Ministry of          Science and Higher Education of the Russian Federation, agreement 075-15-2019-1620 date 08/11/2019 and 075-15-2022-289 date 06/04/2022.

\bibliographystyle{ieeetr}

\bibliography{queue2d}
\end{document}